\documentclass[12pt]{amsart}
\usepackage{amsmath,amssymb,latexsym}
\usepackage{fullpage}
\usepackage{multirow}
\usepackage[colorlinks=true]{hyperref}

\newcommand{\bbZ}{\mathbb{Z}}

\newtheorem{theorem}{Theorem}

\theoremstyle{definition}
\newtheorem*{remark-nonum}{Remark}

\begin{document}

\title[$n$-defective Lehmer pairs\ldots]{$n$-defective Lehmer pairs for small $n$: Corrections and Clarifications}

\author{Paul M. Voutier}
\address{London, UK}
\email{Paul.Voutier@gmail.com}

\date{}

\begin{abstract}
We provide some corrections and clarifications of statements in \cite{Ab} and \cite{BHV}
regarding elements in Lehmer sequences that have no primitive divisors for small $n$.
\end{abstract}


\maketitle

\section{Introduction}

The Primitive Divisor Theorem \cite{BHV} has found many applications in number
theory, especially in the study of diophantine problems. Often knowing that there
is always a primitive divisor of any element $u_{n}$ of a Lucas or Lehmer sequence,
$\left( u_{n} \right)_{n \geq 0}$,
for $n>30$ suffices for applications. But to complete the solution of some
problems, one may need to know what happens for smaller indices, $n$, too.
Those $n$ with only finitely many Lucas ($n>4$ and $n \neq 6$) or Lehmer sequences
($n>6$ and $n \neq 8,10,12$) such that $u_{n}$ does not have a primitive
divisor were treated in \cite{V1}.
The remaining indices $n$, with infinitely many Lucas or Lehmer sequences such
that $u_{n}$ does not have a primitive divisor were treated in Theorem~1.3 of \cite{BHV}.
This was later corrected by Abouzaid \cite{Ab}.

However, in Sept. 2021, the user Seee \cite{Seee} posted a question on mathoverflow about
a Lehmer sequence whose $5$-th element
has no primitive divisor, but does not appear in either \cite[Th\'{e}or\`{e}me~4.1
and the accompanying table on page~313]{Ab} or \cite[Theorem~1.3 and Table~4]{BHV}. In particular,
this user mentioned the Lehmer sequence generated by
$(\alpha, \beta)= \left( \left( 1-\sqrt{5} \right)/2, \left( 1+\sqrt{5} \right)/2 \right)$,

In this note, we show that this sequence was incorrectly omitted from both
\cite{BHV} and \cite{Ab} as well as correcting the previous proofs for $n=5$
and, consequently, also handling correctly $n=10$, adding a missing sequence here.
Furthermore, we revisit the proofs for the other small values of $n$, making
appropriate corrections and changes to the results for such $n$ too. See
Subsection~\ref{subsect:changes} below for a description of all the changes.

The Lucas sequence work for small $n$ in \cite{Ab} has been checked and no errors
were found in Th\'{e}or\`{e}me~4.1 and its associated table there (note that \cite{Ab}
does make some corrections to \cite{BHV} for Lucas sequences for $n=3,4$ and $6$).
So here we focus on the Lehmer sequences.

\subsection{Notation}
\label{subsect:notation}

We start with some notation and facts from \cite{BHV}. Since we  make extensive
use of \cite{BHV} in what follows, we will also indicate when we use equations
from \cite{BHV}. If an equation from \cite{BHV} is used below, we denote it by
(BHV-xyz), where (xyz) is the corresponding equation from \cite{BHV}.
For example, equation~\eqref{eq:bhv-2} below is equation~(2) in \cite{BHV}.

Recall that a \emph{Lehmer pair} is a pair $\left( \alpha, \beta \right)$ of algebraic
integers such that $\left( \alpha+\beta \right)^{2}$ and $\alpha\beta$ are non-zero
coprime rational integers and $\alpha/\beta$ is not a root of unity. For a Lehmer
pair $(\alpha, \beta)$, we define the corresponding sequence of \emph{Lehmer numbers}
by
\begin{equation}
\label{eq:bhv-2}
\tilde{u}_{n} = \tilde{u}_{n}(\alpha, \beta)
= \left\{
\begin{array}{ll}
\dfrac{\alpha^{n}-\beta^{n}}{\alpha-\beta} & \text{if $n$ is odd,} \\[3ex]
\dfrac{\alpha^{n}-\beta^{n}}{\alpha^{2}-\beta^{2}} & \text{if $n$ is even.}
\end{array}
\right.
\tag{BHV-2}
\end{equation}

For a Lehmer pair $(\alpha, \beta)$, a prime number $p$ is a \emph{primitive divisor}
of $\tilde{u}_{n}(\alpha, \beta)$ if $p$ divides $\tilde{u}_{n}$, but does not
divide $\left( \alpha^{2}-\beta^{2} \right)^{2}\tilde{u}_{1} \cdots \tilde{u}_{n-1}$.

A Lehmer pair $(\alpha, \beta)$ such that $\tilde{u}_{n}(\alpha, \beta)$ has no
primitive divisors will be called an \emph{$n$-defective Lehmer pair}.

We say that two Lehmer pairs $\left( \alpha_{1}, \beta_{1} \right)$ and
$\left( \alpha_{2}, \beta_{2} \right)$ are \emph{equivalent} if
$\alpha_{1}/\alpha_{2}=\beta_{1}/\beta_{2} \in \left\{ \pm 1, \pm \sqrt{-1} \right\}$.
A consequence of this that we will use repeatedly is that
\begin{equation}
\label{eq:equiv-cond}
\text{$\left( \alpha_{1}, \beta_{1} \right)$ and
$\left( \alpha_{2}, \beta_{2} \right)$ are equivalent}
\Longrightarrow
\alpha_{1}\beta_{1}= \pm \alpha_{2}\beta_{2}.
\end{equation}
While necessary, this condition
is not sufficient for the pairs to be equivalent.

We will also often use the fact that if two Lehmer pairs $\left( \alpha_{1}, \beta_{1} \right)$
and $\left( \alpha_{2}, \beta_{2} \right)$ are equivalent, then
we have $\tilde{u}_{n}\left( \alpha_{1}, \beta_{1} \right)
= \pm \tilde{u}_{n}\left( \alpha_{2}, \beta_{2} \right)$ for all $n \geq 0$.

As in Section~3 of \cite{BHV}, for a Lehmer pair, $\left( \alpha, \beta \right)$,
we put $p=\left( \alpha+\beta \right)^{2}$, $q=\alpha \beta$ and define $\sqrt{p}$
uniquely as $\alpha+\beta$. So $\alpha$ and $\beta$ are roots of $X^{2}-\sqrt{p}X+q$
and
\[
\alpha, \beta = \frac{\sqrt{p} \pm \sqrt{p-4q}}{2}.
\]

To express our results for $n=5,8,10$ and $12$, we will need to define several
binary recurrence sequences.

For $n=5$ and $10$, we will need the following two sequences.

For $k \geq 0$, let $\phi_{k}$ be the $k$-th element of the Fibonacci sequence
and $\psi_{k}$ be the $k$-th element of the associated sequence, which is sometimes
referred to as the classical Lucas sequence, defined
by $\psi_{0}=2$, $\psi_{1}=1$ and $\psi_{k+1}=\psi_{k}+\psi_{k-1}$.

We will also need to set $\phi_{-2}$, $\psi_{-2}$, $\psi_{-1}$ and $\psi_{-1}$.
Extending the recurrence, $\phi_{k+2}=\phi_{k+1}+\phi_{k}$, so it applies
for all $k \geq -2$, we find that $\phi_{-2}=-1$ and $\phi_{-1}=1$ is consistent
with the other values of $\phi_{k}$.
Similarly, we set $\psi_{-2}=3$ and $\psi_{-1}=-1$.

For $n=8$, we will need the following two sequences.

We define the sequence $\left( \pi_{k} \right)_{k \geq 0}$ by
$\pi_{0}=0$, $\pi_{1}=1$ and $\pi_{k+1}=2\pi_{k}+\pi_{k-1}$.
Similarly, we define the sequence $\left( \rho_{k} \right)_{k \geq 0}$ by
$\rho_{0}=\rho_{1}=1$ and $\rho_{k+1}=2\rho_{k}+\rho_{k-1}$ (the same recurrence
relation as for the $\pi_{k}$'s). To ensure that all the elements used in the
proofs for $n=8$ in both \cite{Ab} and \cite{BHV} are actually defined, we also
set $\rho_{-1}=-1$ and $\pi_{-1}=1$. As with the extension of the sequences for
$n=5$, this extension is consistent with the recurrence relation for these
sequences, along with their $0$-th and $1$-st elements.

For $n=12$, we will need the four sequences, $\left( \zeta_{k}^{(i)} \right)_{k \geq -1}$
for $i=0,1,2,3$. They all satisfy the recurrence
relation $\zeta_{k+1}^{(i)}=4\zeta_{k}^{(i)}-\zeta_{k-1}^{(i)}$ for $i=0,1,2,3$
and have the initial terms in Table~\ref{table:initial-elements}. These are the
same as the sequences
in \cite{BHV}, except that we have included an element for $k=-1$. Similar to
the extensions of the sequences for $n=5$ and $n=10$ above, these additional
elements will ensure that all $12$-defective Lehmer pairs are found.

For those sequences above that we have started at negative indices, we could have
shifted our sequences to have them start at $k=0$ instead, but we have chosen to
start at negative indices so that the $k$-th elements of our sequences equal
the $k$-th elements of the sequences in \cite{BHV}.

\begin{center}
\begin{table}[h]
\begin{tabular}{|c|r|r|r|r|}\hline
$i$ & $0$ & $1$ & $2$ & $3$ \\ \hline
$\zeta_{-1}^{(i)}$ & $-1$ & $2$ & $1$ & $-1$  \\ \hline
 $\zeta_{0}^{(i)}$ &  $0$ & $1$ & $1$ &  $1$  \\ \hline
 $\zeta_{1}^{(i)}$ &  $1$ & $2$ & $3$ &  $5$  \\ \hline
\end{tabular}
\caption{Initial values for $n=12$ sequences}
\label{table:initial-elements}
\end{table}
\end{center}

\subsection{Results}

\begin{theorem}
\label{thm:main}
For $n \in \left\{ 3,4,5,6,8,10,12 \right\}$, up to equivalence, all $n$-defective
Lehmer pairs are of the form
$\left( \left( \sqrt{a}-\sqrt{b} \right)/2, \left( \sqrt{a}+\sqrt{b} \right)/2 \right)$,
where $(a,b)$ is given in Table~$\ref{table:lehmer}$ with $k, \ell, q \in \bbZ$
and $\varepsilon = \pm 1$.
\end{theorem}

\begin{center}
\begin{table}[h]
\begin{tabular}{|l|ll|}\hline
$n$                 & $(a,b)$ & \\ \hline
\multirow{2}{*}{3}  & $(1+q, 1-3q)$, & $q \neq -1,0,1$ \\
                    & $\left( 3^{k}+q, 3^{k}-3q \right)$, & $k>0$, $3 \nmid q$, $(k,q) \neq (1,1)$ \\ \hline
\multirow{2}{*}{4}  & $(1+2q, 1-2q)$, & $q \neq -1,0,1$ \\
                    & $\left( 2^{k}+2q, 2^{k}-2q \right)$, & $k>0$, $2 \nmid q$, $(k,q) \not\in \{ (1,-1), (1,1), (2,1) \}$ \\ \hline
\multirow{2}{*}{5}  & $\left( \phi_{k-2\varepsilon}, \phi_{k-2\varepsilon}-4\phi_{k} \right)$, & $k \geq 3$ \\
                    & $\left( \psi_{k-2\varepsilon}, \psi_{k-2\varepsilon}-4\psi_{k} \right)$, & $k \geq 0$, $(k,\varepsilon) \not\in \left\{ (0,-1), (1,-1) \right\}$ \\ \hline
\multirow{4}{*}{6}  & $(1+3q, 1-q)$, & $q \neq -1,0,1$ \\
                    & $\left( 3^{\ell}+3q, 3^{\ell}-q \right)$, & $\ell>0$, $3 \nmid q$, $(\ell,q) \neq (1,-1)$ \\
                    & $\left( 2^{k}+3q, 2^{k}-q \right)$, & $k>0$, $2 \nmid q$, $(k,q) \neq (1,-1)$ \\
                    & $\left( 2^{k}3^{\ell}+3q, 2^{k}3^{\ell}-q \right)$, & $k,\ell>0$, $\gcd(6,q)=1$ \\ \hline
\multirow{2}{*}{8}  & $\left( \rho_{k-\varepsilon}, \rho_{k-\varepsilon}-4\pi_{k} \right)$, & $k \geq 2$ \\
                    & $\left( 2\pi_{k-\varepsilon}, 2\pi_{k-\varepsilon}-4\rho_{k} \right)$, & $k \geq 2$ \\ \hline
\multirow{2}{*}{10} & $\left( \phi_{k-2\varepsilon}-4\phi_{k}, \phi_{k-2\varepsilon} \right)$, & $k \geq 3$ \\
                    & $\left( \psi_{k-2\varepsilon}-4\psi_{k}, \psi_{k-2\varepsilon} \right)$, & $k \geq 0$, $(k,\varepsilon) \not\in \left\{ (0,-1), (1,-1) \right\}$ \\ \hline
\multirow{3}{*}{12} & $\left( \zeta^{(i)}_{k-\varepsilon}, -\zeta^{(i)}_{k+\varepsilon} \right)$, & $i \in \{ 0,1,2,3 \}$, $k \geq 0$, \\
                    &  & $(i,k, \varepsilon) \not\in \{ (0,0, \pm 1), (0,1, \pm 1), (1,0,\pm 1), (2,0, \pm 1) \}$ \\ \hline
\end{tabular}
\caption{$n$-defective Lehmer pairs}
\label{table:lehmer}
\end{table}
\end{center}

We exclude $n=1$ and $2$ here since $\tilde{u}_{1}=\tilde{u}_{2}=1$, so all
Lehmer pairs are $1$-defective and $2$-defective.

\subsection{Changes from \cite{Ab,BHV}}
\label{subsect:changes}

We provide here the corrections and changes made in this paper.

\subsubsection{$n=3$ changes}

(1) We have excluded $q=-1$ (recall that we put $q=\alpha\beta$ above)
from the pairs $(a,b)=(1+q,1-3q)$. It was included in both \cite{Ab} and
\cite{BHV}. But for $q=-1$, we have $a=0$, so $\left( \alpha+\beta \right)^{2}=0$,
which is not permitted for Lehmer pairs.

\vspace*{1.0mm}

\noindent
(2) We have also corrected a typo in \cite{Ab}. He had
$\left( 3^{k}, 3^{k}-3q \right)$ in his table instead of
$\left( 3^{k}+q, 3^{k}-3q \right)$, although in his proof for $n=3$
this family of pairs is given correctly.

\subsubsection{$n=4$ changes}

(1) We have excluded $q=-1$ from the family of pairs $(a,b)=\left( 1+2q, 1-2q \right)$.
Here we have $(p,q)=(-1,-1)$, so this pair
is excluded by \eqref{eq:bhv-28} below.

\vspace*{1.0mm}

\noindent
(2) We have excluded $(k,q)=(1,-1)$ from the family of pairs $(a,b)=\left( 2^{k}+2q, 2^{k}-2q \right)$.
Here we have $a=0$, so $\left( \alpha+\beta \right)^{2}=0$,
which is not permitted for Lehmer pairs.

Both of these were wrongly included in both \cite{Ab} and \cite{BHV}.

\subsubsection{$n=5$ changes}

(1) We added the pair $(a,b)$ associated with $\left( \psi_{k} \right)_{k \geq -2}$
having $(k,\varepsilon)=(1,1)$. This pair was
omitted from both \cite{Ab} and \cite{BHV} -- this is the example found by
Seee \cite{Seee}.

\vspace*{1.0mm}

\noindent
(2) Also for the pairs associated with $\left( \psi_{k} \right)_{k \geq -2}$,
$(k,\varepsilon)=(0,-1)$
must be excluded since the resulting Lehmer pair, $(a,b)=\left( \psi_{2}, \psi_{2}-4\psi_{0} \right)
=(3,-5)$, is equal to the Lehmer pair for $(k,\varepsilon)=(0,1)$, where
$(a,b)=\left( \psi_{-2}, \psi_{-2}-4\psi_{0} \right)=(3,-5)$.

\subsubsection{$n=6$ changes}

(1) We have excluded $q=-1$ from the family of pairs $(a,b)=\left( 1+3q, 1-q \right)$.
This was included in both \cite{Ab} and \cite{BHV}. When $q=-1$, $p=1+3q=-2$ and
the pair $(p,q)=(-2,-1)$ is excluded by \eqref{eq:bhv-28} below.

\vspace*{1.0mm}

\noindent
(2) $(k,q)=(1,-1)$ for the family of pairs $\left( 2^{k}+3q, 2^{k}-q \right)$
and $(\ell,q)=(1,-1)$ for the family of pairs $\left( 3^{\ell}+3q, 3^{\ell}-q \right)$ are
excluded both here and in \cite{Ab}, but both were included incorrectly in
\cite{BHV}.

\subsubsection{$n=8$ changes}

For $n=8$, there are no changes from \cite{Ab} and \cite{BHV}.

\subsubsection{$n=10$ changes}

(1) The pair $(a,b)$ associated with the
sequence $\left( \psi_{k} \right)_{k \geq -2}$ having $(k,\varepsilon)=(1,1)$
is included here. It was omitted from both \cite{Ab} and \cite{BHV}.

\vspace*{1.0mm}

\noindent
(2) For the pairs associated with $\left( \psi_{k} \right)_{k \geq -2}$,
$(k,\varepsilon)=(0,-1)$ must be excluded. It was wrongly included in both
\cite{Ab} and \cite{BHV}.

\subsubsection{$n=12$ changes}

(1) There are six sequences defined in \cite{Ab}. There are two sequences
associated with each of what we call $\left( \zeta_{k}^{(2)} \right)$ and
$\left( \zeta_{k}^{(3)} \right)$. There is one sequence for $\varepsilon=-1$
and another for $\varepsilon=1$. Labelling them here as $\left( \zeta_{k}^{(i,-1)} \right)$
and $\left( \zeta_{k}^{(i,1)} \right)$ for $i=2,3$, we have the relationships
$\zeta_{k}^{(i,1)}=-\zeta_{k+1}^{(i,-1)}$ for $i=2,3$.
Hence the pairs in \cite{Ab} (see his table on page~313) arising from $i \in \{ 2,3\}$,
$k \geq 1$ and $\varepsilon=-1$ are equivalent to those arising from $i$, $k-1$
and $\varepsilon=1$. So here we exclude the sequences (and pairs) in \cite{Ab}
arising from his $\varepsilon=-1$ to ensure that none of the Lehmer
pairs arising from the pairs $(a,b)$ for $n=12$ in our theorem are equivalent
to other pairs there. This also keeps the four sequences here for $n=12$ consistent
with the four sequences used for $n=12$ in \cite{BHV}.

\vspace*{1.0mm}

\noindent
(2) For $n=12$ in Table~4 on page~79 of \cite{BHV}, the excluded values
should be as follows.

For $\varepsilon=1$, exclude
$k=0,1$ for all $i$ (since $k-2\varepsilon<0$ and these sequences are not defined
at negative indices in \cite{BHV}), and also $(i,k)=(0,2)$ (since $a=0$ in this case).

For $\varepsilon=-1$, exclude $(i,k)=(0,0)$ (since $b=0$ in this case).

\vspace*{1.0mm}

\noindent
(3) Also for $n=12$ in Table~4 of \cite{BHV}, the pair $(a,b)=(1,5)$ is missing.
This arises from $i=3$ and $k=1$ with $\varepsilon=-1$ in the notation of \cite{Ab}.
We include it here by extending the sequence of $\zeta_{k}^{(3)}$'s to include
$\zeta^{(3)}_{-1}=-1$.

\vspace*{1.0mm}

\noindent
(4) In the proof of Theorem~1.3 for $n=12$ in \cite{BHV}, the correct values of
$p$ and $q$ are obtained and stated on line~-5 of page~90. However, when determining
$a$ and $b$ on lines~-2 and -1 of page~90, it appears that $(a,b)=(q-4p,q)$ (rather
than $(a,b)=(p,p-4q)$) is used to obtain the pair in their Table~4. The proof on
page~90 of \cite{BHV} with the values of $p$ and $q$ there actually gives the pair
$(a,b)=(p,p-4q)=\left( \zeta^{(i)}_{k-\varepsilon}, -\zeta^{(i)}_{k+\varepsilon} \right)$
given in the entry for $n=12$ in the table on page~313 of \cite{Ab}. These are
the pairs we provide in our result.
However, since
both collections of pairs $(a,b)$ give rise to the same collection of Lehmer pairs.

\subsection{Code}

PARI/GP \cite{Pari} code for work associated with this paper can be found in the
\verb!primitive-divisors! subdirectory of the
\url{https://github.com/PV-314/sequences} repository.
For the Lehmer sequences, the code finds all pairs $(\alpha, \beta)$ obtained
from $(a,b)$ with $|a| \leq 10^{7}$ for $n=5,8,10$ and $12$.
The author found this code
useful as a check that no further pairs were overlooked as well as for realising
that for $n=12$, the four sequences (extended as above) in \cite{BHV} sufficed.

Similarly, there is PARI/GP code for checking the results for Lucas sequences for
$n=2,3,4$ and $6$.

The author is
very happy to help interested readers who have any questions, problems or
suggestions for the use of this code.

\section{Proof of Theorem~\ref{thm:main}}

%
We have
\begin{equation}
\label{eq:bhv-26}
q \neq 0,
\tag{BHV-26}
\end{equation}
\begin{equation}
\label{eq:bhv-27}
\gcd \left( \tilde{u}_{n}, q \right) = \gcd \left( \Phi_{n}(p,q), q \right)
= \gcd (p,q) =1,
\tag{BHV-27}
\end{equation}
\begin{equation}
\label{eq:bhv-28}
(p,q) \not\in \left\{ \pm (1,1), \pm (2,1), \pm (3,1), \pm (4,1) \right\}.
\tag{BHV-28}
\end{equation}

\subsection{$n=3$}

The arguments in both \cite{Ab} and \cite{BHV} correctly yield
$(a,b)=(1+q,1-3q)$ with $q \in \bbZ$; or else
$(a,b)=\left( 3^{k}+q, 3^{k}-3q \right)$ with $k \geq 1$ and $q \in \bbZ$.

By \eqref{eq:bhv-26}, we exclude $q=0$ from both families.

\subsubsection{$(a,b)=(1+q,1-3q)$}

Consider the Lehmer pair $\left( \alpha, \beta \right)$ obtained from
$(a,b)=(1+q,1-3q)$. If $q=-1$, then $a=0$, which we also exclude. If $q=1$, then
$p=2$. This pair $(p,q)$ is excluded by \eqref{eq:bhv-28}.
Since $p=1+q$, we have $\gcd(p,q)=1$ for other values of $q$.

For a Lehmer pair $\left( \alpha, \beta \right)$ obtained from $(a,b)=(1+q,1-3q)$,
we have $\alpha\beta=(a-b)/4=((1+q)-(1-3q))/4=q$. So by \eqref{eq:equiv-cond}, two pairs
$\left( \alpha_{1}, \beta_{1} \right)$
and $\left( \alpha_{2}, \beta_{2} \right)$ associated with distinct $\left( 1+q_{1}, 1-3q_{1} \right)$
and $\left( 1+q_{2}, 1-3q_{2} \right)$ can only be equivalent if $q_{2}=-q_{1}$.

If the two Lehmer pairs are equivalent, then $\tilde{u}_{4}\left( \alpha_{1}, \beta_{1} \right)=1-q_{1}
=\pm \tilde{u}_{4}\left( \alpha_{2}, \beta_{2} \right)=\pm \left( 1-q_{2} \right)$.
Since $q_{2}=-q_{1}$, we must have $1-q_{1}=\pm \left( 1+ q_{1} \right)$. The
only solution is $q_{1}=0$, which we have already excluded. So none of the Lehmer
pairs obtained from any $(a,b)=(1+q,1-3q)$ is equivalent with another such pair.

\subsubsection{$(a,b)=\left( 3^{k}+q, 3^{k}-3q \right)$}

Consider the Lehmer pairs $\left( \alpha, \beta \right)$ obtained from
$(a,b)=\left( 3^{k}+q, 3^{k}-3q \right)$, where $k \geq 1$ and $q \in \bbZ$.
Here $p=3^{k}+q$ and so $\gcd(p,q)|3^{k}$.
This equals $1$, as required, only if $q \not\equiv 0 \pmod{3}$.

If $q=1$, then \eqref{eq:bhv-28} excludes $k=1$.

As with the first family (i.e., those where $(a,b)=(1+q,1-3q)$), here two pairs
$\left( \alpha_{1}, \beta_{1} \right)$ and $\left( \alpha_{2}, \beta_{2} \right)$
obtained from the distinct pairs
$\left( 3^{k_{1}}+q_{1}, 3^{k_{1}}-3q_{1} \right)$ and
$\left( 3^{k_{2}}+q_{2}, 3^{k_{2}}-3q_{2} \right)$ can only be equivalent if
$q_{2}=-q_{1}$. As above, we look at the $4$-th element to show that such Lehmer
pairs are not equivalent. We can assume that $1 \leq k_{1}<k_{2}$. We require
$\tilde{u}_{4}\left( \alpha_{1}, \beta_{1} \right)=3^{k_{1}}-q_{1}
=\pm \tilde{u}_{4}\left( \alpha_{2}, \beta_{2} \right)=\pm \left( 3^{k_{2}}-q_{2} \right)
=\pm \left( 3^{k_{2}}+q_{1} \right)$.

From $3^{k_{1}}-q_{1}=3^{k_{2}}+q_{1}$, we obtain $3^{k_{1}} \left( 3^{k_{2}-k_{1}} - 1 \right)
=-2q_{1}$. Since $k_{1}>0$, this implies that $3|q_{1}$, which we exclude.

From $3^{k_{1}}-q_{1}=- \left( 3^{k_{2}}+q_{1} \right)$, we obtain
$3^{k_{1}}=3^{k_{2}}$. But this implies $k_{1}=k_{2}$.

So none of the Lehmer pairs associated with any $\left( 3^{k}+q, 3^{k}-3q \right)$
is equivalent with another such pair.

Lastly, we note that $\tilde{u}_{3}=3^{k}>1$ for such Lehmer pairs. So none of
the Lehmer pairs associated with any $\left( 3^{k}+q, 3^{k}-3q \right)$ is
equivalent with a Lehmer pair associated with any $\left( 1+q, 1-3q \right)$,
since we have $\tilde{u}_{3}=1$ for them.

\subsection{$n=4$}

The arguments in both \cite{Ab} and \cite{BHV} both correctly yield the two families
given. The only mistake is not excluding $(k,q)=(1,-1)$ (where $(a,b)=(0,4)$)
from the second family, since $p=a=0$ does not give rise to a Lehmer pair.

The proof that none of these pairs is equivalent to another is the same as for
$n=3$.

\subsection{$n=5$}

The argument for $n=5$ in \cite{BHV} is correct up to and including equations
(34) and (35) there, which we repeat there:
\begin{align}
\label{eq:bhv-34}
p = & \varepsilon^{k+1} \frac{\eta^{\varepsilon k-2}-\overline{\eta}^{\varepsilon k-2}}{\sqrt{5}}
= \frac{\eta^{k-2\varepsilon}-\overline{\eta}^{k-2\varepsilon}}{\sqrt{5}}
= \phi_{k-2\varepsilon}, \tag{BHV-34} \\
q = & \varepsilon^{k+1} \frac{\eta^{\varepsilon k}-\overline{\eta}^{\varepsilon k}}{\sqrt{5}}
= \frac{\eta^{k}-\overline{\eta}^{k}}{\sqrt{5}}
= \phi_{k}, \nonumber
\end{align}
or
\begin{align}
\label{eq:bhv-35}
p = & \varepsilon^{k} \left( \eta^{\varepsilon k-2}+\overline{\eta}^{\varepsilon k-2} \right)
= \eta^{k-2\varepsilon}+\overline{\eta}^{k-2\varepsilon}
= \psi_{k-2\varepsilon}, \tag{BHV-35} \\
q = & \varepsilon^{k} \left( \eta^{\varepsilon k}+\overline{\eta}^{\varepsilon k} \right)
= \eta^{k}+\overline{\eta}^{k}
= \psi_{k}. \nonumber
\end{align}

Recall from page~88 immediately after equation~(32) in \cite{BHV} that
$\varepsilon=\pm 1$ and $k$ is a non-negative integer.

For $k \geq 3$, we have $\phi_{k} \geq 2$. In the case of equation~\eqref{eq:bhv-34},
we have $q=\phi_{k}$, so only $k=0,1$ and $2$ possibly lead to violations of
\eqref{eq:bhv-28}. Also, $q=\phi_{0}=0$, so $k=0$ is excluded by \eqref{eq:bhv-26}.
In the case of equation~\eqref{eq:bhv-35}, among the non-negative values of $k$,
only $k=1$ possibly leads to violations of \eqref{eq:bhv-28}.

We now consider these values of $k$ in each case.

For \eqref{eq:bhv-34} with $k=1$, we have $k-2\varepsilon=-1,3$. So $q=\phi_{1}=1$,
while $p=\phi_{-1}=1$ or $p=\phi_{3}=2$. Both of these are excluded by
\eqref{eq:bhv-28}. With $k=2$, we have $k-2\varepsilon=0,4$. So $q=\phi_{2}=1$,
while $p=\phi_{0}=0$ or $p=\phi_{4}=3$. We do not allow $p=0$ (recall from our
definition of Lehmer pairs that $p=\left( \alpha+\beta \right)^{2} \neq 0$) and
$(p,q)=(3,1)$ is excluded by \eqref{eq:bhv-28}. Hence for equation~\eqref{eq:bhv-34},
we need only consider $k \geq 3$.

For \eqref{eq:bhv-35} with $k=1$, we have $k-2\varepsilon=-1,3$. So $q=\psi_{1}=1$,
while $p=\psi_{-1}=-1$ or $p=\psi_{3}=4$. The second of these is excluded by
\eqref{eq:bhv-28}. However, $(p,q)=(-1,1)$ is permissible. This solution gives
rise to the Lehmer pair $\left( \alpha, \beta \right) = \left( \left( \sqrt{-1} + \sqrt{-5} \right)/2,
\left( \sqrt{-1} - \sqrt{-5} \right)/2 \right)$. We find that
$\tilde{u}_{n}(\alpha,\beta)=0,1,1,-2,-3,5$ for $n=0,1,2,3,4,5$.
Since $\left( \alpha^{2}-\beta^{2} \right)^{2}=5$, by the definition of a
primitive divisor, we see that $\tilde{u}_{5}(\alpha,\beta)$ does not have a
primitive divisor.

From $a=p$ and $b=p-4q$, we obtain the pairs $(a,b)$ in the theorem for $n=5$.

From equations~\eqref{eq:bhv-26}--\eqref{eq:bhv-28}, the conditions required for
$(\alpha, \beta)$ to be a Lehmer pair hold.

We now show that if
$\left( \alpha_{1}, \beta_{1} \right)$ and $\left( \alpha_{2}, \beta_{2} \right)$
are distinct pairs from the statement of the theorem, then they are not equivalent.
Let $\left( p_{1}, q_{1} \right)$
and $\left( p_{2}, q_{2} \right)$ be the pairs $(p,q)$ associated with
$\left( \alpha_{1}, \beta_{1} \right)$ and $\left( \alpha_{2}, \beta_{2} \right)$
respectively.
If they are equivalent, then $\alpha_{1}\beta_{1}=\left( a_{1}-b_{1} \right)/4
= \pm \alpha_{2}\beta_{2}= \pm \left( a_{2}-b_{2} \right)/4$.

If $\left( a_{i}, b_{i} \right) = \left( \phi_{k_{i}-2\varepsilon}, \phi_{k_{i}-2\varepsilon}-4\phi_{k} \right)$,
for some non-negative integer $k_{i}$, then $\alpha_{i}\beta_{i}=\phi_{k_{i}}$.

Similarly, if
$\left( a_{i}, b_{i} \right) = \left( \psi_{k_{i}-2\varepsilon}, \psi_{k_{i}-2\varepsilon}-4\psi_{k} \right)$,
for some non-negative integer $k_{i}$, then $\alpha_{i}\beta_{i}=\psi_{k_{i}}$.

So if $\left( \alpha_{1}, \beta_{1} \right)$ and $\left( \alpha_{2}, \beta_{2} \right)$
are equivalent, then there exist $k_{1}, k_{2} \geq 0$ such that
$\phi_{k_{1}}=\pm \phi_{k_{2}}$ for $k_{1} \neq k_{2}$,
or $\psi_{k_{1}}=\pm \psi_{k_{2}}$ for $k_{1} \neq k_{2}$
or $\phi_{k_{1}}=\pm \psi_{k_{2}}$.

If $\phi_{k_{1}}=\pm \phi_{k_{2}}$ for $k_{1} \neq k_{2}$, then it must be the
case that one
of $k_{1}$ and $k_{2}$ is $1$, while the other is $2$, since the sequence of
$\phi_{k}$'s is strictly increasing for $k \geq 2$.
In either case, $q_{1}=q_{2}=1$,
while $p_{1}, p_{2} \in \left\{ \phi_{-1}=1, \phi_{0}=0, \phi_{3}=2, \phi_{4}=3 \right\}$.
Both $p_{1}$ and $p_{2}$ must be non-zero. Equation~\eqref{eq:bhv-28} eliminates
the other possibilities here.

If $\psi_{k_{1}}=\pm \psi_{k_{2}}$ for $k_{1} \neq k_{2}$, then we proceed in
the same way.

If $\phi_{k_{1}}=\pm \psi_{k_{2}}$, then $k_{1} \in \{ 0,1,2,3,4 \}$ (i.e.,
$q_{1}=\phi_{k_{1}} \in \{ 0,1,2,3\}$) and $k_{2} \in \{ 0,1,2 \}$ (i.e.,
$q_{2}=\psi_{k_{2}} \in \{ 1,2,3 \}$).

But since $q_{1} \neq 0$ by \eqref{eq:bhv-26}, we have
$k_{1} \in \{ 1,2,3,4 \}$
$q_{1}=\phi_{k_{1}} \in \{ 1,2,3\}$. Checking each of the possibilities for $p_{1}$
and $q_{1}$, we are left with $\left( p_{1}, q_{1} \right) \in \left\{ (1,2), (5,2),
(1,3), (8,3) \right\}$.

$q_{1}=q_{2}=2$ implies that $p_{2}=3$, while $p_{1}=1$ or $5$.

$q_{1}=q_{2}=3$ implies that $p_{2}=2$ or $7$, while $p_{1}=1$ or $8$.

Checking each of these possibilities, we find that none of the pairs are equivalent.
This completes the proof for $n=5$.


\subsection{$n=6$}

There are no changes to the proofs in \cite{Ab,BHV} except for the observations for $n=6$ in Subsection~\ref{subsect:changes}.

\subsection{$n=8$}

There are no changes to the proofs in \cite{Ab,BHV} except to check $k=0$
and $\varepsilon=1$ in both families of pairs, $(a,b)$, in Table~\ref{table:lehmer}.
For the first family, we
have $p=\rho_{-1}=-1$ and $q=\pi_{0}=0$. This is excluded by \eqref{eq:bhv-26}.
For the second family, $p=2\pi_{-1}=2$ and $q=\rho_{0}=1$. This is excluded by
\eqref{eq:bhv-28}.

On page~90 of \cite{BHV} (last paragraph before $n=10$ case), it is stated that
\eqref{eq:bhv-28} shows that $k \geq 2$.
However for $k=1$ and $\varepsilon=1$, we need more than \eqref{eq:bhv-28}.
Here we have either $p=\rho_{k-\varepsilon}=\rho_{0}=1$
and $q=\pi_{k}=\pi_{0}=0$, which we exclude by \eqref{eq:bhv-26}, or
$p=2\pi_{k-\varepsilon}=2\pi_{0}=0$ and $q=\rho_{k}=\rho_{1}=1$, which we exclude
since $p=0$ does not yield a Lucas pair.

\subsection{$n=10$}

As stated in the treatment of $n=10$ on page~90 of \cite{BHV}, $\Phi_{10}(\alpha,\beta)
=\Phi_{5}(-\alpha, \beta)$. Hence from the expressions for $\alpha$ and $\beta$
in terms of $a$ and $b$, a pair $(a,b)$ appears for $n=10$ if and only if $(b,a)$
appears for $n=5$.

\subsection{$n=12$}

For $n=12$, the analysis in \cite{Ab} is correct. However, there are equivalent
Lehmer pairs among the pairs provided in \cite{Ab}. In \cite{BHV}, the analysis
is also correct up to, and including line~-3 on page~90, although we must extend
the sequences, as we have done here, so that $p=\zeta_{k-\varepsilon}^{(i)}$ is
defined when $k=0$ and $\varepsilon=1$.

When $k=0$ and $\varepsilon=1$, we have
$(p,q)=\left( \zeta_{k-\varepsilon}^{(i)}, \zeta_{k}^{(i)} \right)
=\left( \zeta_{-1}^{(i)}, \zeta_{0}^{(i)} \right)=(-1,0),(2,1)$ and
$(1,1)$ for $i=0,1,2$, respectively. By \eqref{eq:bhv-26}, we cannot have $q=0$.
By \eqref{eq:bhv-28}, the pairs, $(p,q)$, for $i=1$ and $2$ are eliminated too.

When $k=0$ and $\varepsilon=-1$, we have
$(p,q)=\left( \zeta_{1}^{(i)}, \zeta_{0}^{(i)} \right)=(1,0),(2,1)$ and
$(3,1)$ for $i=0,1,2$, respectively. Again, all of these are eliminated by
\eqref{eq:bhv-26} and \eqref{eq:bhv-28}.

For $i=0$, $k=1$ and $\varepsilon=1$, we have
$(p,q)=\left( \zeta_{0}^{(0)}, \zeta_{1}^{(0)} \right)=(0,1)$. This is
eliminated because $p=(\alpha+\beta)^{2}$ must be non-zero.

For $i=0$, $k=1$ and $\varepsilon=-1$, we have
$(p,q)=\left( \zeta_{2}^{(0)}, \zeta_{1}^{(0)} \right)=(4,1)$. This is
eliminated by \eqref{eq:bhv-28}.

We now consider equivalent Lehmer pairs.

Recall from \eqref{eq:equiv-cond} that a necessary condition for the Lehmer
pairs $\left( \alpha_{1}, \beta_{1} \right)$ and $\left( \alpha_{2}, \beta_{2} \right)$
to be equivalent is that $\alpha_{1}\beta_{1}=\pm \alpha_{2}\beta_{2}$.
Furthermore, if $\alpha_{j}=\left( \sqrt{a_{j}}+\sqrt{b_{j}} \right)/2$ and
$\beta_{j}=\left( \sqrt{a_{j}}-\sqrt{b_{j}} \right)/2$, then
$\alpha_{j}\beta_{i}=\left( a_{j}-b_{j} \right)/4$.
Here this quantity is $\left( \zeta^{(i)}_{k+1}+\zeta^{(i)}_{k-1} \right)/4
=\zeta^{(i)}_{k}$, by
the recurrence relation. Since for each $i$, the elements, $\zeta^{(i)}_{k}$,
are distinct for $k \geq 1$, there are no equivalent Lehmer pairs arising within
each sequence, $\left( \zeta^{(i)}_{k} \right)$ for distinct indices, $k_{1}, k_{2} \geq 1$.

It remains to show that there are no equivalent pairs across sequences. We showed
above that if a Lehmer pair, $(\alpha,\beta)$, comes from $\zeta^{(i)}$, then
$\alpha\beta=\zeta^{(i)}_{k}$ for some $k \geq 0$. We can only have equivalent
Lehmer pairs if there are values in one such sequence that are $\pm$ values in
another such sequence. The elements of these sequences with non-negative
indices are all non-negative, so we need only consider when $\zeta^{(i_{1})}_{k_{1}}
=\zeta^{(i_{2})}_{k_{2}}$. The only value that occurs in multiple such sequences is $1$, which occurs as
$\zeta^{(0)}_{1}, \zeta^{(1)}_{0}, \zeta^{(2)}_{0}$ and $\zeta^{(3)}_{0}$.
This can be seen by observing
that for $k \geq 1$,
\[
\zeta_{k-1}^{(3)}<\zeta_{k}^{(1)}<\zeta_{k}^{(2)}<\zeta_{k+1}^{(0)}<\zeta_{k}^{(3)},
\]
a fact that follows by an inductive argument.

But we exclude the pairs associated with the first three of $\zeta^{(0)}_{1}$,
$\zeta^{(1)}_{0}, \zeta^{(2)}_{0}$ and $\zeta^{(3)}_{0}$.
Hence, there are no such equivalent pairs.

\end{document}